# On regularity of solutions to Poisson's equation

Rahul Garg
rgarg@tx.technion.ac.il

Daniel Spector
dspector@tx.technion.ac.il

Department of Mathematics
Technion - Israel Institute of Technology

**Abstract**

In this note, we announce new regularity results for some locally integrable distributional solutions to Poisson's equation. This includes, for example, the standard solutions obtained by convolution with the fundamental solution. In particular, our results show that there is no qualitative difference in the regularity of these solutions in the plane and in higher dimensions.

## 1 Main Results

In this note, we announce new regularity results for certain $L^1_{loc}(\mathbb{R}^N)$ solutions to Poisson's equation

$$-\Delta u = f. \tag{1.1}$$

To be more precise, let $f \in L^p(\mathbb{R}^N)$ for some $\frac{N}{2} < p \leq N$, so that we know there exists $u \in L^1_{loc}(\mathbb{R}^N)$ which satisfies (1.1) in the sense of distributions (see, for example, [2][Chapter 6, Theorem 6.21, p.157]). This means that

$$-\int_{\mathbb{R}^N} u \Delta \varphi = \int_{\mathbb{R}^N} f \varphi, \tag{1.2}$$

for all $\varphi \in C_c^\infty(\mathbb{R}^N)$. If we restrict ourselves to solutions which satisfy the growth condition

$$\frac{u(x)}{|x|} \to 0 \text{ as } |x| \to \infty, \tag{1.3}$$

then, up to a constant, such a $u \in L^1_{loc}(\mathbb{R}^N)$ is unique. In fact, when $supp\ f$ is compact, this solution can simply be written as convolution with the fundamental solution to Laplace's equation (the logarithmic potential for $N = 2$ and the Newtonian potential for $N \geq 3$). We, however, do not impose such an assumption on $f$ in this article.

The main result we announce here is the following theorem on the regularity of this solution.

**Theorem 1.1** *Suppose $N \geq 2$ and $\frac{N}{2} < p \leq N$. If $\frac{N}{2} < p < N$, assume $f \in L^p(\mathbb{R}^N)$, while if $p = N$ assume $f \in L^N(\mathbb{R}^N) \cap L^q(\mathbb{R}^N)$ for some $1 < q < N$. Then there exists $u \in L^1_{loc}(\mathbb{R}^N)$ which satisfies (1.2), for which one has the following regularity estimates.*

i) If $\frac{N}{2} < p < N$, then
$$|u(x) - u(z)| \leq C|x-z|^{2-\frac{N}{p}} \|f\|_{L^p(\mathbb{R}^N)}.$$

ii) If $p = N$, then
$$|u(x) - u(z)| \leq C|x-z| (|\ln|x-z|| + 1)^{\frac{1}{N'}} \left(\|f\|_{L^N(\mathbb{R}^N)} + \|f\|_{L^q(\mathbb{R}^N)}\right).$$

*In particular, this solution satisfies* (1.3)*, and hence it is unique up to a constant.*

An important feature of our result is that the estimates we obtain are uniform and with sharp exponent. We here explain how such a result cannot be obtained as not a consequence of the known embeddings. For this, first recall that one has the following inclusions
$$W^{2,p}(\mathbb{R}^N) \subsetneq \{u \in L^1_{loc}(\mathbb{R}^N) : \Delta u \in L^p(\mathbb{R}^N)\} \subsetneq W^{2,p}_{loc}(\mathbb{R}^N).$$

While it is known that the uniform estimates we obtain hold for functions in $W^{2,p}(\mathbb{R}^N)$, for which the exponents are known to be sharp (see, for example, [3][Chapter 1, p. 62, Remark 1] and [1][Corollary 5]), in general the solution need not be in $W^{2,p}(\mathbb{R}^N)$. For example, taking $f = \chi_{B(0,1)}(x)$, one can verify that when $N = 2$ we have $u \notin L^q(\mathbb{R}^2)$ for any $1 \leq q \leq +\infty$, while when $N = 3$, $u \notin L^q(\mathbb{R}^3)$ for $1 \leq q \leq 3$. Further, the above inclusion implies that the solution is in the space $W^{2,p}_{loc}(\mathbb{R}^N)$, so that one could apply standard embeddings to obtain local versions of our results. However, such embeddings do not give one uniform estimates, where the constant $C$ is independent of $x, z \in \mathbb{R}^N$. A simple example of this is $x^2 \in W^{2,p}_{loc}(\mathbb{R}^N)$.

Another interesting aspect of our result is that in the case $N = 2$, it is based on a new representation formula for the solution. Precisely, if we define the map
$$\tilde{T}_j h(x) := \frac{1}{2\pi} \int_{\mathbb{R}^2} \left[\frac{y_j - x_j}{|y-x|} - \frac{y_j}{|y|}\right] h(y) \, dy,$$

then we claim
$$u := \tilde{T}_1 R_1 f + \tilde{T}_2 R_2 f \tag{1.4}$$

solves (1.2), where $R_j$ is the standard $j$-th Riesz transform,
$$R_j f(x) := \frac{1}{2\pi} \int_{\mathbb{R}^2} \frac{x_j - y_j}{|x-y|^3} f(y) \, dy.$$

The regularity and uniqueness of $u$ are then a consequence of the following theorem on the mapping properties of $\tilde{T}_j$.

**Theorem 1.2** *Let $1 < p \leq 2$.*

i) *If $1 < p < 2$, then there exists $C = C(p)$ such that*
$$|\tilde{T}_j h(x) - \tilde{T}_j h(z)| \leq C|x-z|^{2-\frac{2}{p}} \|h\|_{L^p(\mathbb{R}^2)}$$

*for all $h \in L^p(\mathbb{R}^2)$ and $j = 1, 2$.*

*ii)* If $p = 2$ and $1 \leq q < 2$, then there exists $C = C(q)$ such that

$$|\tilde{T}_j h(x) - \tilde{T}_j h(z)| \leq C|x-z|(|\ln|x-z|| + 1)^{\frac{1}{2}} \left(\|h\|_{L^2(\mathbb{R}^2)} + \|h\|_{L^q(\mathbb{R}^2)}\right)$$

for all $h \in L^2(\mathbb{R}^2) \cap L^q(\mathbb{R}^2)$ and $j = 1, 2$.

A third remark concerning Theorem 1.1 is that while when $N \geq 3$ the regularity results for the case $\frac{N}{2} < p < N$ are already known (see [4][Section 4.2, Theorem 2.2, p. 155]), the case $p = N$ has not previously been treated in this setting. Here, as we previously alluded to, the correct embedding for Sobolev functions had been understood by Brezis and Wainger [1][Corollary 5]. As we are not, in general, under the same hypothesis, the regularity result from Theorem 1.1 must be deduced otherwise. We therefore require the following theorem concerning the mapping properties of the modified Newtonian potential.

**Theorem 1.3** *For any $1 \leq q < N$, there exists $C = C(q, N)$ such that*

$$|\tilde{I}_2 f(x) - \tilde{I}_2 f(z)| \leq C|x-z|(|\ln|x-z|| + 1)^{\frac{1}{N'}} \left(\|f\|_{L^q(\mathbb{R}^N)} + \|f\|_{L^N(\mathbb{R}^N)}\right)$$

*for all $f \in L^N(\mathbb{R}^N) \cap L^q(\mathbb{R}^N)$.*

Here, we have defined the modified Newtonian potential

$$\tilde{I}_2 f(x) := \frac{\Gamma(\frac{N}{2})}{2\pi^{\frac{N}{2}}(N-2)} \int_{\mathbb{R}^N} \left[\frac{1}{|y-x|^{N-2}} - \frac{1}{|y|^{N-2}}\right] f(y) \, dy,$$

since the Newtonian potential need not be well-defined on $L^p(\mathbb{R}^N)$ for $\frac{N}{2} \leq p \leq N$.

We now sketch a proof that the function defined by (1.4) solves (1.2). First, we remark that for $\frac{N}{2} < p < N$, one can show the estimate

$$\int_{\mathbb{R}^2} \left|\frac{y_j - x_j}{|y-x|} - \frac{y_j}{|y|}\right| |R_j f(y)| \, dy \leq C|x|^{2-\frac{N}{p}} \|f\|_{L^p(\mathbb{R}^2)},$$

for $j = 1, 2$. Therefore, Fubini's theorem implies

$$-\int_{\mathbb{R}^2} u \Delta\varphi = -\frac{1}{2\pi} \sum_{j=1,2} \int_{\mathbb{R}^2} \left(\int_{\mathbb{R}^2} \left[\frac{y_j - x_j}{|y-x|} - \frac{y_j}{|y|}\right] R_j f(y) \, dy\right) \Delta\varphi(x) \, dx$$

$$= -\frac{1}{2\pi} \sum_{j=1,2} \int_{\mathbb{R}^2} \left(\int_{\mathbb{R}^2} \left[\frac{y_j - x_j}{|y-x|} - \frac{y_j}{|y|}\right] \Delta\varphi(x) \, dx\right) R_j f(y) \, dy.$$

Further, since the divergence theorem implies $\int_{\mathbb{R}^2} \Delta\varphi(x) \, dx = 0$, we have that

$$\int_{\mathbb{R}^2} \left[\frac{y_j - x_j}{|y-x|} - \frac{y_j}{|y|}\right] \Delta\varphi(x) \, dx = \int_{\mathbb{R}^2} \frac{y_j - x_j}{|y-x|} \Delta\varphi(x) \, dx.$$

Now, we define

$$g_j(y) := -\frac{1}{2\pi} \int_{\mathbb{R}^2} \frac{y_j - x_j}{|y-x|} \Delta\varphi(x) \, dx.$$

If we can show that $g_j = R_j\varphi$ almost everywhere, then we would have

$$-\int_{\mathbb{R}^2} u\Delta\varphi = \sum_{j=1,2} \int_{\mathbb{R}^2} R_j\varphi R_j f$$
$$= \int_{\mathbb{R}^2} f\varphi,$$

which is the thesis. Notice that

$$g_j(y) = -\frac{y_j}{2\pi}\int_{\mathbb{R}^2} \frac{1}{|y-x|}\Delta\varphi(x)\,dx + \frac{1}{2\pi}\int_{\mathbb{R}^2} \frac{1}{|y-x|}x_j\Delta\varphi(x)\,dx,$$

and therefore,

$$\widehat{g_j}(\xi) = \frac{1}{2\pi i}\frac{\partial}{\partial \xi_j}\left((2\pi|\xi|)^{-1}\widehat{\Delta\varphi}(\xi)\right) + (2\pi|\xi|)^{-1}\left(\widehat{x_j\Delta\varphi}(\xi)\right)$$
$$= i\left[\frac{\partial}{\partial \xi_j}(|\xi|\widehat{\varphi}(\xi)) - \frac{1}{|\xi|}\frac{\partial}{\partial \xi_j}(|\xi|^2\widehat{\varphi}(\xi))\right]$$
$$= -i\frac{\xi_j}{|\xi|}\widehat{\varphi}(\xi)$$
$$= \widehat{R_j\varphi}.$$

Here, the above should be interpreted in the sense of tempered distributions, and with the convention

$$\widehat{\varphi}(\xi) = \int_{\mathbb{R}^2} \varphi(x)e^{-2\pi i x\cdot\xi}\,dx$$

for the Fourier transform. Thus, we have proved that $g_j = R_j\varphi$ as distributions, which implies almost everywhere equality as functions, and the result is demonstrated.

Finally, we mention that the proofs of Theorems 1.2 and 1.3 will appear in a forthcoming work, where we also address regularity properties of more general cases of Riesz and Riesz-type potentials, as well as the application of these results to deduce the embedding theorem of Brezis and Wainger [1][Corollary 5] in the supercritical case.

## Acknowledgements

The authors would like to thank Yehuda Pinchover and Georgios Psaradakis for their helpful discussions during the preparation of this work, as well as Igor Verbitsky for valuable comments on an early draft of the manuscript. The first author is supported in part by a research grant (No: 471/13) of Amos Nevo from the Israel Science Foundation and a postdoctoral fellowship from the Planning and Budgeting Committee of the Council for Higher Education of Israel. The second author is supported in part by a Technion Fellowship.

[2] E.H. Lieb, M. Loss, Analysis, Second edition, American Mathematical Society, Providence, RI, 2001.

[3] V. Maz'ya, Sobolev spaces, Springer-Verlag, Berlin, 1985.

[4] Y. Mizuta, Potential Theory in Euclidean Spaces, Gakkōtosho, Tokyo, 1996.

[5] E. Stein, Singular Integrals and Differentiability Properties of Functions, Princeton University Press, New Jersey, 1970.5